\documentclass[12pt,intlim,righttag]{article}
\usepackage{amsmath,amsthm}
\usepackage{amssymb}
\usepackage{amsfonts}

\newcommand{\la}{\langle}
\newcommand{\ra}{\rangle}

\newcommand{\td}{\tilde}

\newcommand{\be}{\begin}
\newcommand{\ee}{\end}
\newcommand{\lbl}{\label}

\newcommand\beq{\begin{equation}}
\newcommand\eeq{\end{equation}}

\newcommand{\beaa}{\begin{eqnarray*}}
\newcommand{\eeaa}{\end{eqnarray*}}

\theoremstyle{Theorem}

\theoremstyle{corollary}

\theoremstyle{remark}

\theoremstyle{definition}

\def\dbF{\mathbb{F}}

\def\e{\varepsilon}

\def\f{\varphi}

\def\o{\omega}

\def\dbF{\mathbb{F}}

\def\f{\varphi}

\def\L{\Lambda}
\def\O{\Omega}

\begin{document}
\title{On Change of Variable Formulas for non-anticipative functionals
 }

\author{M. Mania$^{1)}$ and R. Tevzadze$^{2)}$}

\date{~}
\maketitle

\begin{center}
$^{1)}$ A. Razmadze Mathematical Institute of Tbilisi State University, 6 Tamarashvili Str., Tbilisi 0177; and
Georgian-American University, 8 Aleksidze Str., Tbilisi 0193, Georgia,
\newline(e-mail: mania@rmi.ge)
\\
$^{2)}$ Georgian-American University, 8 Aleksidze Str., Tbilisi 0193, Georgia,
Georgian Technical Univercity, 77 Kostava str., 0175,
Institute of Cybernetics,  5 Euli str., 0186, Tbilisi,
Georgia
\newline(e-mail: rtevzadze@gmail.com)
\end{center}

\begin{abstract}
{\bf Abstract.}
For non-anticipative functionals, differentiable in Chitashvili's sense, the It\^o formula
for cadlag semimartingales is proved. Relations between different notions of functional derivatives
are established.
\end{abstract}

\bigskip


\


\section{Introduction}

    The classical It\^o \cite{ito} formula shows that for a sufficiently smooth function\\ $(f(t,x), t\ge0, x\in R)$ the transformed process
$f(t,X_t)$
is a semimartingale for any semimartingale $X$ and provides a decomposition of the process $f(t,X_t)$ as a sum of stochastic
 integral relative to $X$ and a process of finite variation. This formula is applicable to functions of the current value of semimartingales,
but in many applications, such as statistics of random processes, stochastic optimal control or mathematical finance, uncertainty affects
through the whole history of the process and it is necessary to consider functionals of entire path of a semimartingale.

In 2009 Dupire \cite{Dupire} proposed a method to extend the It\^o formula for non-anticipative functionals
 using naturally defined pathwise  time and space derivatives. The space derivative measures the sensitivity of a functional $f:D([0,T], R)\to R$
to a variation in the endpoint of a path $\o\in D([0,T], R)$ and is defined as a limit
$$
\partial_\o f(t,\o)=\lim_{h\to 0}\frac{f(t,\o+hI_{[t,T]})-f(t,\o)}{h},
$$
if this limit exists, where $D([0,T])$ is the space of RCLL ( right continuous with left limits) functions. Similarly is defined the second order space
derivative $\partial_{\omega\omega}f:= \partial_{\omega}(f_{\omega}).$

The definition of the  time derivative is based on the flat extension of a path $\o$  up to time $t+h$ and is defined as a  limit
$$
\partial_t f(t,\o)=\lim_{h\to 0+}\frac{f(t+h,\o^t)-f(t,\o)}{h},
$$
whenever this limit exists, where $\o^t=\o(.\wedge t)$ is the path of $\o$ stopped at time $t$.

 If a continuous non-anticipative  functional $f$ is from $C^{1,2}$ , i.e.,  if $\partial_t f, \partial_\o f$, $\partial_{\o\o}f$ exist and are continuous
 with respect to the metric $d_\infty$ (defined in section 2)
 and $X$ is a continuous semimartingale, Dupire \cite{Dupire} proved that the process $f(t,X)$ is also a semimartingale and
$$
f(t, X)=f(0, X)+\int_0^t\partial_t f(s,X)ds+\int_0^t\partial_\o f(s,X)dX_s
$$
\begin{equation}\label{itoc}
+\frac{1}{2}\int_0^t\partial_{\o\o}f(s, X)d\la X\ra_s.
\end{equation}
For the special case of $f(t,X_t)$ these derivatives coincide with the usual space and time derivatives and the above formula
 reduces to the standard It\^o formula.
Erlier related works are the works by Ahn  \cite{ahn} and
Tevzadze \cite{T2}, where It\^o's formula was derived in very particular cases of functionals that assume the knowledge of the whole path
without path dependent dynamics.
Further works extending this theory and  corresponding references one can see   in \cite{CF1}, \cite{CF2}, \cite{LScS},\cite{O}.

Motivated by applications in stochastic optimal control,  before Dupire's work, Chitashvili (1983) defined differentiability of non-anticipative
functionals in a different way and proved the corresponding It\^o  formula for continuous semimartingales. His definition is based
 on "hypothetical" change of variable formula for continuous functions of finite variation.

We formulate Chitashvili's definition of differentiability and  present his change of variable formula in a simplified form and  for one-dimensional case.

Let $C_{[0,T]}$ be the space of continuous functions on $[0,T]$ equipped with the uniform norm. Let $f(t,\o)$ be non-anticipative continuous mapping of $C_{[0,T]}$ into $C_{[0,T]}$ and
denote by ${\cal V}_{[0,T]}$ the space of functions of finite variation on $[0,T]$.

A continuous  non-anticipative functional  $f$ is differentiable  if  there exist  continuous functionals $f^0$ and $f^1$ such that  for all  $\o\in C_{[0,T]}\cap {\cal V}_{[0,T]}$
\begin{equation}\label{chd}
f(t,\o)=f(0,\o)+\int_0^tf^0(s,\o)ds+\int_0^tf^1(s,\o)d\o_s.
\end{equation}
A functional $f$ is two times differentiable if $f^1$ is differentiable, i.e., if
there exist  continuous functionals $f^{0,1}$ and $f^{1,1}$ satisfying
\begin{equation}\label{chd2}
f^1(t,\o)=f^1(0,\o)+\int_0^tf^{1,0}(s,\o)ds+\int_0^tf^{1,1}(s,\o)d\o_s.
\end{equation}
 for all  $\o\in C_{[0,T]}\cap {\cal V}_{[0,T]}$.

Here functionals $f^0, f^1$ and $f^{1,1}$ play the role of time, space and the second order space derivatives respectively.

It was proved by Chitashvili \cite{Ch} that if the functional $f$ is two times differentiable then the process $f(t,X)$ is a semimartingale for any continuous semimartingale $X$ and is represented as
$$
f(t, X)=f(0, X)+\int_0^tf^0(s,X)ds+\int_0^tf^1(s,X)dX_s
$$
\begin{equation}\label{itoc}
+\frac{1}{2}\int_0^tf^{1,1}(s, X)d\la X\ra_s.
\end{equation}

The idea of the proof of change of variable formula (\ref{itoc}) for semimartingales is to use  the change of variable formula for functions of finite variations,
first for the function $f$ and then for its derivative $f^1$, before approximating a continuous semimartingale $X$ by processes of finite variation.

In the paper Ren et al \cite{RTZ} a wider class of $C^{1,2}$  functionals  was proposed, which is
based on the Ito formula itself. We formulate this definition in equivalent form and  in one-dimensional case.

The function $f$ belongs to $C^{1,2}_{RTZ}$, if $f$  is a continuous non-anticipative functional on $[0,T]\times C_{[0,T]}$ and there exist  continuous non-anticipative functionals
$\alpha, z, \gamma$, such that
\begin{equation}\label{itoc2}
f(t, X)=f(0, X)+\int_0^t\alpha(s,X)ds+\int_0^tz(s,X)dX_s +\frac{1}{2}\int_0^t\gamma(s, X)d\la X\ra_s
\end{equation}
for any continuous semimartingale $X$.
The functionals $\alpha, z$ and $\gamma$ also play  the role of time, first and second order space derivatives  respectively.

Since any process of finite variation is a semimartingale and any deterministic semimartingale is a function of finite variation, it follows from $f\in C^{1,2}_{RTZ}$ that
 $f$ is differentiable in the Chitashvili sense and
\begin{equation}\label{ChT}
\alpha=f^0,\;\;\;z=f^1.
\end{equation}
Becides, any $C^{1,2}$ process in the Dupire or Chitashvili  sense is in $C^{1,2}_{RTZ}$,
which is a consequence of the functional It\^o formula proved in    \cite{Dupire} and   \cite{Ch} respectively.
Although, the definition of the class $C^{1,2}_{RTZ}$ does not require that $\gamma$ be (in some sense) the derivative of $z$, but
 if $f\in C^{1,2}$  in the Chitashvili sense,
then beside equality  (\ref{ChT}) we also have that $\gamma=f^{1,1}$ (i.e., $\gamma=z^1$).

 Our goal is to extend the formula (\ref{itoc}) for RCLL (or cadlag in French terminology)  semimartingales and to establish how Dupire's, Chitashvili's and other derivatives are related.

Since the bumped path used in the definition of Dupire's vertical derivative is not continuous even if $\o$ is continuous,
to compare derivatives defined by (\ref{chd}) with Dupire's derivatives, one should extend Chitashvili's definition to RCLL processes, or
to modify Dupire's derivative in such a way that perturbation of continuous paths remain continuous.

The direct extension of Chitashvili's definition of differentiability for  RCLL functions is  following:

A  continuous  functional  $f$ is differentiable, if  there exist  continuous  functionals $f^0$ and $f^1$ (continuous with respect to the metric
$d_\infty$ defined by (\ref{rho}))
such that  $ f(\cdot,\o)\in {\cal V}_{[0,T]}$ for all  $\o\in  {\cal V}_{[0,T]}$ and
\begin{equation}\label{xvii}
f(t,\o)=f(0,\o)+\int_0^tf^0(s,\o)ds+\int_0^tf^1(s-,\o)d\o_s
\end{equation}
$$
+\sum_{s\le t}\big[f(s,\o)-f(s-,\o)-f^1(s-,\o)\Delta\o_s\big],
$$
for all  $(t,\o)\in [0,T]\times {\cal V}_{[0,T]}$.

In order to compare Dupire's derivatives with Chitashvili's  derivatives, we introduce another type of
 vertical  derivative where, unlike to Dupire's derivative $\partial_\omega f$,  the path deformation of
continuous paths  are also  continuous.

We say that a  non-anticipative functional $f(t,\omega)$ is vertically differentiable   and denote this differential by $D_\omega f(t,\omega)$, if
the limit
\begin{equation}
D_\omega f(t,\omega):=\lim_{h\to0, h>0}\frac{f(t+h,\o^{t}+\chi_{t,h})-f(t+h,\o^{t})}{h},
\end{equation}
exists  for all $(t,\o)\in [0,T]\times {D}_{[0,T]}$, where
$$
\chi_{t,h}(s)=(s-t)1_{(t,t+h]}(s)+h1_{(t+h,T]}(s).
$$
Let  $f(t,\o)$  be  differentiable in the sense of (\ref{xvii}).
Then, as proved in Proposition 1,
\begin{equation}
f^0(t,\o)=\partial_t f(t,\omega)\;\;\;\;\text{and}\;\;\;\; f^1(t,\o)=D_\omega f(t,\omega).
\end{equation}
 for all  $(t,\o)\in [0,T]\times {D}_{[0,T]}$.

Thus,  $f^0$ coincides with Dupire's time derivative, but $f^1$  is equal to $D_\o f$ which is different from Dupire's vertical derivative in general. The simplest counterexample
is $f(t,\o)=\o_t-\o_{t-}$. It is evident that  in this case $\partial_\o f=1$ and $D_\o f=0$. In general, if $g(t,\o):=f(t-,\o)$ then
$D_\o g(t,\o)=D_\o f(t,\o)$ and $\partial_\o g(t,\o)=0$ if corresponding derivatives of $f$ exist. However, under stronger conditions, e.g.  if $f\in C^{1,1}$ in the Dupire  sense, then $D_\o f$ exists and $D_\o f=f^1=\partial_\o f.$

The paper is organized as follows: In section 2 we extend Citashvili's change of variable formula for RCLL semimartingales and
give an application of this formula on the convergence of ordinary integrals to the stochastic integrals. In section 3 we establish relations between different type
of derivatives for non-anticipative functionals.

\section{The It\^o formula  according to Chitashvili for cadlag semimartingales}

Let $\O:= D([0,T], R)$ be the set of c\`{a}dl\`{a}g
paths. Denote by $\omega$  the elements of $\O$, by $\o_t$ the value of $\o$ at time $t$ and let  $\o^t=\o(\cdot\wedge t)$ be the path of $\o$ stopped at  $t$.
Let  $B$ be the canonical process defined by $B_t(\o)=\o_t$,   $\dbF=(F_t,t\in[0,T])$ the  corresponding filtration  and  let
$\L:= [0,T]\times\O$.

The functional $f:[0,T]\times D[0,T]\to R$ is non-anticipative if
$$
f(t,\o)=f(t,\o^t)
$$
for all $\o\in D[0,T]$, i.e., the process $f(t,\o)$ depends only on the path of $\o$ up to time $t$ and is $\dbF$- adapted.

Following Dupire, we define semi-norms on $\O$ and a pseudo-metric on $\L$ as follows:
for any $(t, \omega), ( t', \omega') \in\L$,
%
\begin{eqnarray}
\label{rho}
\|\omega\|_{t}&:=& \sup_{0\le s\le t} |
\omega_s|,\nonumber\\[-8pt]\\[-8pt]
d_\infty\bigl((t, \omega),\bigl(
t', \omega'\bigr) \bigr)&:=& \bigl|t-t'\bigr| +
\sup_{0\le s\le T} \bigl|\omega_{t\wedge s} - \omega'_{t'\wedge
s}\bigr|.\nonumber
\end{eqnarray}
Then $(\O, \|\cdot\|_{T})$ is a Banach space and $(\L, d_\infty)$ is
a complete pseudo-metric space.
Let ${\cal V}={\cal V}[0,T]$ be the set of finite variation paths from $\O$. Note that,
if $f\in C(\L)$, then from $\Delta \o_t=0$ follows $f(t,\o)-f(t-,\o)=0$, since $d_\infty((t_n,\o),(t,\o))\to 0$ when $t_n\uparrow t$.
Hence $f(t,\o)-f(t-,\o)\neq 0$ means $\Delta \o_t\neq 0$.

Note that any functional $f:[0,T]\times\Omega\to R$ continuous with respect to $d_\infty$ is non-anticipative. In this paper we consider only $d_\infty$-continuous, and hence non-anticipative,
functionals.

{\bf {Definition 1.}}
We say that a continuous functional  $f\in C([0,T]\times \O)$ is differentiable , if  there exist $f^0\in C([0,T]\times \O)$ and $f^1\in C([0,T]\times \O)$
such that  for all  $\o\in {\cal V}$ the process $ f(t,\o)$ is of finite variation and
\begin{equation}\label{xv}
f(t,\o)=f(0,\o)+\int_0^tf^0(s,\o)ds+\int_0^tf^1(s-,\o)d\o_s
\end{equation}
$$
+\sum_{s\le t}\big[f(s,\o)-f(s-,\o)-f^1(s-,\o)\Delta\o_s\big],
$$
for all  $(t,\o)\in [0,T]\times\cal V$.

A functional $f$ is two times differentiable if $f^1$ is differentiable, i.e., if
there exist $f^{0,1}\in C([0,T]\times \O)$ and $f^{1,1}\in C([0,T]\times \O)$
such that  for all  $(t,\o)\in [0,T]\times\cal V$
\begin{equation}\label{two}
f^1(t,\o)=f^1(0,\o)+\int_0^tf^{1,0}(s,\o)ds+\int_0^tf^{1,1}(s-,\o)d\o_s + V^1(t,\o),
\end{equation}
where
$$
V^1(t,\o)=\sum_{s\le t}\big(f^1(s,\o)-f^1(s-,\o)-f^{1, 1}(s-,\o)\Delta\o_s\big).
$$

Now we give a generalization of Theorem 2  from Chitashvili \cite{Ch}
  for general  cadlag (RCLL) semimartingales.

\be{thr}
 Let $f$ be two times differentiable in the sense of Definition 1 and assume that for some $K>0$
\begin{equation}\label{v}
|f(t,\o)-f(t-,\o)-f^1(t-,\o)\Delta\o_t|\le K(\Delta\o_t)^2,\;\; \forall\o\in\cal V.
\end{equation}
 Then for any semimartingale $X$ the process $f(t,X)$ is a semimartingale and
$$
f(t, X)=f(0, X)+\int_0^tf^0(s,X)ds+\int_0^tf^1(s-,X)dX_s
$$
\begin{equation}\label{ito}
+\frac{1}{2}\int_0^tf^{1,1}(s, X)d\la X^c\ra_s+\sum_{s\le t}\big[f(s,X)-f(s-,X)-f^1(s-,X)\Delta X_s\big].
\end{equation}
\ee{thr}
{\it Proof.} Let first assume that $X$ is  a semimartingale with the decomposition
\begin{equation}\label{dec0}
X_t=A_t+M_t, t\in[0,T],
\end{equation}
where $M$ is a continuous local martingale and $A$ is a process of finite variation having only  finite number of jumps, i.e., the jumps of $A$ are exhausted by
graphs of finite number of stopping times $(\tau_i, 1\le i\le l, l<\infty)$.

Let $X_t^n= A_t+M^n_t$ and
\begin{equation}
M^n_t= n\int_0^tM_s\exp(-n(\la M\ra_t-\la M\ra_s)d\la M\ra_s.
\end{equation}
It is proved in  \cite{Ch} that
\begin{equation}\label{mc}
\sup_{s\le t}|M^n_s-M_t|\to 0, \;\;\;as\;\;n\to\infty,\;\;\; a.s.
\end{equation}

Since $X^n$ is of bounded variation, $f$ is differentiable   and $\Delta X^n_t=\Delta A_t=\Delta X_t$, it follows from (\ref{xv}) that
$$
f(t,X ^n)=f(0, X)+\int_0^tf^0(s,X^n)ds
$$
$$
+\int_0^tf^1(s-,X^n)dX_s +\int_0^tf^1(s-,X^n)d(M^n_s-M_s)
$$
\begin{equation}\label{itod}
+\sum_{s\le t}\big(f(s,X^n)-f(s-,X^n)-f^1(s-,X^n)\Delta X_s\big).
\end{equation}

Since $X$ admits finite number of jumps, by continuity of $f$ and $f^1$,

\begin{equation}\label{jumpb}
\sum_{s\le t}\big(f(s,X^n)-f(s-,X^n)-f^1(s-,X^n)\Delta X_s\big)\to
\end{equation}
$$
\to\sum_{s\le t}\big(f(s,X)-f(s-,X)-f^1(s-,X)\Delta X_s\big)
$$

The continuity  of $f, f^0, f^1$ and  relation (\ref{mc}) imly that
\begin{equation}\label{fxn1a}
f(t,X^n)\to f(t,X),\;\;\;as\;\;n\to\infty,\;\;\; a.s.,
\end{equation}
\begin{equation}\label{fxn22a}
\int_0^tf^0(s,X^n)ds\to\int_0^tf^0(s,X)ds\;\;\;as\;\;n\to\infty,\;\;\; a.s..
\end{equation}
by the dominated convergence theorem and
\begin{equation}\label{fxn4a}
\int_0^tf^1(s-,X^n)dX_s\to\int_0^tf^1(s-,X)dX_s\;\;\;as\;\;n\to\infty,\;\;\; a.s..
\end{equation}
by the dominated convergence theorem for stochastic integrals.
Here we may use the dominated convergence theorem, since by continuity of $f^i ( i=0,1)$ the process
$\sup_{n, s\le t}|f^i(s-, X^n)|$ is locally bounded (see Lemma A1).

Let us show now that
\begin{equation}\label{fx12aa}
\int_0^tf^1(s-,X^n)d(M^n_s-M_s)\to\frac{1}{2}\int_0^tf^{1,1}(s,X)d\la M\ra_s.
\end{equation}
Integration by parts and (\ref{two}) give
$$
\int_0^tf^1(s,X^n)d(M^n_s-M_s)=(M^n_t-M_t)f^1(t,X^n)-
$$
$$
-\int_0^t(M^n_s-M_s)f^{1,0}(s,X^n)ds-\int_0^t(M^n_s-M_s)f^{1,1}(s-,X^n)dA_s
$$
$$
-\int_0^t(M^n_s-M_s)f^{1,1}(s-,X^n)dX^n_s-\int_0^t(M^n_s-M^c_s)dV^1(s, X^n)=
$$
\begin{equation}\label{i5}
=I^1_t(n)+I^2_t(n)+I^3_t(n)+I^4_t(n) +I_t^5(n).
\end{equation}

$I^1_t(n)\to 0$ (as $n\to\infty$, a.s.) by continuity of $f^1$ and (\ref{mc}).

$I^2_t(n)$ and  $I_t^3(n)$ tend to zero (as $n\to\infty$, a.s.) by continuity of $f^{1,0}$  and $f^{1,1}$,  relation (\ref{mc}) and  by the dominated convergence theorem
(using the same  arguments as in (\ref{fxn22a})-(\ref{fxn4a})).

Moreover, since  $A$ admits finite number of jumps at  $(\tau_i, 1\le i\le l)$
\begin{equation}\label{jump2}
I_t(5)=\sum_{s\le t}(M^n_s-M_s)\big(f^1(s,X^n)-f^1(s-,X^n)-f^{1,1}(s-,X^n)\Delta A_s\big)
\end{equation}
$$
=\sum_{i\le l}(M^n_{\tau_i}-M_{{\tau_i}})\big(f^1(\tau_i,X^n)-f^1(\tau_i-,X^n)-f^{1,1}(\tau_i-,X^n)\Delta A_{\tau_i}\big)
$$
$$
\le \sup_{s\le t}|M^n_s-M_s|\big(2l\sup_{n, s\le t}|f^1(s,X^n)|+\sup_{n, s\le t}|f^{1,1}(s,X^n)|\sum_{i\le l}|\Delta A_{\tau_i}|\big)\to 0,
$$
as $n\to\infty$, since the continuity of $f^1, f^{1,1}$, relation (\ref{mc}) and Lemma A1 imply that
 $\sup_{n, s\le t}|f^1(s,X^n)|+\sup_{n, s\le t}|f^{1,1}(s,X^n)|<\infty$ (a.s.)

Let us consider now the term
$$
I_t^4(n)=\int_0^t(M_s-M^n_s)f^{1,1}(s,X^n)dM^n_s
$$
Let
$$
K^n_t=\int_0^t(M_s-M^n_s)dM^n_s.
$$
Using the formula of integration by parts we have
$$
K^n_t=-\frac{1}{2}(M_t^n)^2+M_t M_t^n-\int_0^tM_s^ndM_s
$$
and it follows from  (\ref{mc}), the dominated convergence theorem and equality $M_t^2=2\int_0^tM_sdM_s+\la M\ra_t$,  that
\begin{equation}\label{kn}
sup_{s\le t}|K^n_s-\frac{1}{2}\la M\ra_s|\to 0, \;\;\;as\;\;n\to\infty,\;\;\; a.s.
\end{equation}
From definition of $M^n$, using the formula of integration by parts, it  follows that $M^n$ admits  representation
$$
M^n_t=n\int_0^t(M_s-M^n_s)d\la M\ra_s.
$$
Therefore
$$
K^n_t=n\int_0^t(M_s-M^n_s)^2d\la M\ra_s.
$$
This  implies that $K^n$ is a sequence of increasing processes, which is stochastically bounded  by (\ref{kn}) (i.e. satisfies the condition UT from (\cite{JMP})
and by  theorem 6.2 of(\cite{JMP}) (it follows also from lemma 12 of \cite{CF1})
$$
\int_0^t(M_s-M^n_s)f^{1,1}(s,X^n)dM^n_s=
$$
$$
=\int_0^tf^{1,1}(s,X^n)dK^n_s\to\frac{1}{2}\int_0^tf^{1,1}(s,X)d\la M\ra_s,\;\;\;n\to\infty,
$$
which (together with (\ref{i5}))  implies the convergence (\ref{fx12aa}).
Therefore,  the formula  (\ref{ito})
for the process $X$ with decomposition (\ref{dec0}) follows by passage to the limit in (\ref{itod}) using relations (\ref{jumpb})-(\ref{fx12aa}).
Note that in this cased the condition (\ref{v}) is not needed.

Let consider now the general case.  Any semimartingale $X$ admits a decomposition     $X_t=A_t+M_t$, where $A$ is a process of finite variation and $M$ is a locally square integrable martingale
(such decomposition is not unique, but the continuous martingale parts coincide for all such decompositions of $X$, which is sufficient for our goals) see \cite{J}.
Let $M_t=M_t^c+M^d_t$, where $M^c$ and $M^d$ are continuous and purely discontinuous martingale parts of $M$
respectively. Let $A_t=A_t^c+A_t^d$ be the decomposition of $A$, where $A^c$ and $A^d$ are continuous  and purely discontinuous processes  of finite variations respectively.
Note that $A^d$ is the sum of its jumps, whereas $M^d$ is the sum of compensated jumps of $M$. So,
we shall use the decomposition
\begin{equation}\label{dec1}
X_t=A_t^c+A_t^d+M_t^c+M_t^d
\end{equation}
for $X$ and using localization arguments, without loss of generality, one can assume that $M^c$ and $M^d$ are square integrable martingales.

Let $M^d_t(n)$ be the compensated sum of jumps of $M$ of amplitude greater than $1/n$, which is a martingale of finite variation and
is expressed as a difference
\begin{equation}\label{jump}
M^d_t(n)=B^n_t-\widetilde{ B_t^n},
\end{equation}
where $B^n_t=\sum_{s\le t}\Delta M_sI_{(|\Delta M_s|\ge 1/n)}$ and $\widetilde{B^n}$ is the dual predictable projection of $B^n$. It
can be expressed also  as compensated stochastic integral
(see \cite{DM})
$$
M^d_t(n)=\int_0^tI_{(|\Delta M_s|>\frac{1}{n})}{}_{\overset{\bullet}C}dM_s,
$$
where by $H{}_{\overset{\bullet}C}Y$ we denote the compensated stochastic integral.
Since
$$
M^d_t(n)-M_t^d=\int_0^tI_{(0<|\Delta M_s|\le\frac{1}{n})}{}_{\overset{\bullet}C}dM_s,
$$
it follows from Doob's inequality and from \cite{DM} (theorem 33, Ch.VIII) that
$$
E\sup_{s\le t}|M_s^d(n)-M_s^d|^2\le const E[M^d(n)-M^d]_t= const E[I_{(0<|\Delta M|\le\frac{1}{n})}{}_{\overset{\bullet}C}M]
$$
$$
\le const E\int_0^tI_{(0<|\Delta M_s|\le\frac{1}{n})}d[M]_s\to 0, \;\;\;as\;\;n\to\infty
$$
by dominated convergence theorem, since $E[M^d]_T<\infty$.
Hence
\begin{equation}\label{md}
\sup_{s\le t}|M^d_s(n)-M^d_s|\to 0, \;\;\;as\;\;n\to\infty,\;\;\; a.s.
\end{equation}
for some subsequence, for which we preserve the same notation.

Let
$$
A_t^d(n)=\sum_{s\le t}I_{(|\Delta A_s|>\frac{1}{n})}\Delta A_s=\int_0^tI_{(|\Delta A_s|>\frac{1}{n})}dA_s.
$$
Since
$$
|A^d_t-A_t^d(n)|\le\int_0^tI_{(0<|\Delta A_s|\le\frac{1}{n})}|dA_s|
$$
we have that
\begin{equation}\label{ad}
\sup_{s\le t}|A^d_s(n)-A^d_t|\to 0, \;\;\;as\;\;n\to\infty,\;\;\; a.s.
\end{equation}

Let
$$
X^n_t= A^c_t+A_t^d(n)+M_t^d(n)+M_t^c.
$$
Relations (\ref{md}) and (\ref{ad}) imply that
\begin{equation}\label{x}
\sup_{s\le t}|X_s(n)-X_s|\to 0, \;\;\;as\;\;n\to\infty,\;\;\; a.s.,
\end{equation}
Thus,  $X^n$ is a sum of continuous local martingale $M^c$ and a process of finite variation $A^c_t+A_t^d(n)+M_t^d(n)$
which admits only finite number of jumps for every $n\ge 1$.

Therefore, as it is already proved,
$$
f(t,X^n)=f(0,X^n)+\int_0^tf^0(s,X^n)ds+\int_0^tf^1(s-,X^n)dX_s
$$
$$
+\int_0^tf^1(s-,X^n)d(M_s^n(d)-M_s^d)+\int_0^tf^1(s-,X^n)d(A_s^n(d)-A_s^d)
$$
$$
+\frac{1}{2}\int_0^tf^{1,1}(s, X)d\la X^c\ra_s
$$
\begin{equation}\label{fxnv}
+\sum_{s\le t}\big(f(s,X^n)-f(s-,X^n)-f^1(s-,X^n)\Delta X^n_s\big).
\end{equation}

By continuity  of $f, f^0$ and $f^1$
\begin{equation}\label{fxn1}
f(t,X^n)\to f(t,X),\;\;\;as\;\;n\to\infty,\;\;\; a.s.,
\end{equation}
\begin{equation}\label{fxn22}
\int_0^tf^0(s,X^n)ds\to\int_0^tf^0(s,X)ds\;\;\;as\;\;n\to\infty,\;\;\; a.s..
\end{equation}
by the dominated convergence theorem and
\begin{equation}\label{fxn4}
\int_0^tf^1(s-,X^n)dX_s\to\int_0^tf^1(s-,X)dX_s\;\;\;as\;\;n\to\infty,\;\;\; a.s..
\end{equation}
by the dominated convergence theorem for stochastic integrals (using the same arguments as in (\ref{fxn22a})- (\ref{fxn4a})).

By properties of compensated stochastic integrals
$$
\int_0^tf^1(s-,X^n)d(M^d_s(n)-M^d_s)=\int_0^tf^1(s-,X^n)I_{(0<|\Delta M_s|\le\frac{1}{n})}{}_{\overset{\bullet}C}dM_s
$$
and using theorem 33, Ch. VIII from \cite{DM}
$$
E\big(\int_0^tf^1(s-,X^n)I_{(0<|\Delta M_s|\le\frac{1}{n})}{}_{\overset{\bullet}C}dM_s\big)^2
$$
\begin{equation}\label{fx}
\le const E\int_0^t(f^1(s-,X^n))^2I_{(0<|\Delta M_s|\le\frac{1}{n})}d[M^d]_s\to 0\;\;\;as\;\;n\to\infty
\end{equation}
by dominated convergence theorem, since
 $\sup_{n, s\le t}(f^1(s,X^n))^2$ is locally bounded (by Lemma A1 from appendix), $I_{(0<|\Delta M_s|\le\frac{1}{n})}\to 0$ and $E[M^d]_T<\infty$.

Similarly,  $\int_0^tf^1(s-,X^n)d(A_s^n(d)-A_s^d)$ also tends to zero, since
\begin{equation}\label{fxan}
\int_0^tf^1(s-,X^n)d(A_s^n(d)-A_s^d)\le \int_0^t|f^1(s-,X^n)|I_{(0<|\Delta A_s|\le\frac{1}{n})}|dA_s|\to 0.
\end{equation}

From (\ref{jump})
$$
\Delta M^n_s(d)=\Delta M_sI_{(|\Delta M_s|\ge 1/n)} - \big( \Delta MI_{(|\Delta M|\ge 1/n)}\big)_s^p,
$$
where $Y^p$ is the usual projection of $Y$. Here we used the fact that the jump of the dual projection of $B^n$ is
the usual  projection of the jump, i.e. $\Delta\widetilde{B^n_t}=(\Delta B^n)_t^p$. Therefore,
using condition (\ref{v})  we have that
$$
|(f(s,X^n)-f(s-,X^n)-f^1(s-,X^n)\Delta X^n_s|\le const. (\Delta X^n_s)^2
$$
$$
= const. \big(\Delta A_sI_{(|\Delta A_s|\ge 1/n)}+\Delta M_sI_{(|\Delta M_s|\ge 1/n)} - ( \Delta MI_{(|\Delta M|\ge 1/n)})_s^p\big)^2
$$
\begin{equation}\label{jump2}
\le 3 const.\big( (\Delta A_s)^2+(\Delta M_s)^2+ E( (\Delta M_s)^2/F_{s-})\big).
\end{equation}

Since, it follows from (\ref{x}) and continuity of $f$ and $f^1$, that
$$
f(s,X^n)-f(s-,X^n)-f^1(s-,X^n)\Delta X^n_s\to f(s,X)-f(s-,X)-f^1(s-,X)\Delta X_s
$$
and
$$
\sum_{s\le t}\big ( (\Delta A_s)^2+(\Delta M_s)^2+ E( (\Delta M_s)^2/F_{s-})\big ) < \infty,
$$
the dominated convergence theorem implies that
$$
\sum_{s\le t}\big(f(s,X^n)-f(s-,X^n)-f^1(s-,X^n)\Delta X^n_s\big)
$$
\begin{equation}\label{jump3}
\to\sum_{s\le t}\big(f(s,X)-f(s-,X)-f^1(s-,X)\Delta X_s\big), \;\;\;as\;\;\;n\to\infty.
\end{equation}

Therefore, passing to the limit in (\ref{fxnv}) it follows from (\ref{fxn1})-(\ref{jump3}) that (\ref{ito}) holds.\qed

Now we give one application of the change of variable formula (\ref{ito}) to the convergence of stochastic integrals.
If $g(t,x), t\ge0, x\in R)$ is a function of two variables admitting continuous partial derivatives $\partial g(t,x)/\partial t$, $\partial g(t,x)/\partial x$
and $V^n$ is a sequence of processes of finite variations
converging to the Wiener process, then it was proved by Wong and Zakai \cite{WZ} that the sequence of ordinary integrals $\int_0^tg(s,V^n_s)dV^n_s$
converges to the Stratanovich stochastic integral. The following assertion generalizes this result for non-anticipative functionals $g(t,\o)$.

{\bf Corollary}. Assume that $f(t,\o)$ is differentiable in the sense of Definition 1 and there is a continuous on $[0,T]\times D([0,T])$ functional $F(t,\o)$ such that
\begin{equation}\label{str0}
F(t,\o)=\int_0^tf(s-,\o)d\o_s
\end{equation}
For all $\o\in{\cal V}_{[0,T]}$.
Let $X$ be a cadlag semimartingale and let $(V^n,n\ge1)$ be a sequence of processes of finite variation converging to $X$ uniformly on $[0, T]$.
Then
\begin{equation}\label{str}
\lim_{n\to\infty}\int_0^tf(s-, V^n)dV^n_s= \int_0^tf(s-,X)dX_s
+\frac{1}{2}\int_0^tf^{1}(s, X)d\la X^c\ra_s.
\end{equation}
Proof:  By continuity of $F$ and (\ref{str0})
\begin{equation}\label{str1}
\lim_{n\to\infty}\int_0^tf(s-, V^n)dV^n_s=\lim_{n\to\infty}F(t,V^n)=F(t,X).
\end{equation}
It is evident that
$$
F^1(t,\o)=f(t,\o),\;\; F^0(t,\o)=0\;\;\text{and}\;\;\;F(t,\o)-F(t-,\o)-F^1(t-,\o)\Delta\o_t=0,
$$
Thus,  $F$ is two times differentiable in the sense of definition 1 and condition (\ref{v}) is automatically satisfied.
Therefore, by the It\^o formula (\ref{ito})
$$
 F(t,X)=\int_0^tf(s-,X)dX_s
+\frac{1}{2}\int_0^tf^{1}(s, X)d\la X^c\ra_s,
$$
which, together with (\ref{str1}) implies the convergence (\ref{str}).

\section{The relations between various definitions of functional derivatives}

Following Dupire \cite{Dupire} we define time and space derivatives, called also horizontal and vertical derivatives of the non-anticipative functionals.

{\bf Definition 2}.
A  non-anticipative functional $f(t,\omega)$ is said to be horizontally differentiable at $(t,\omega)\in\Lambda$ if the limit
\begin{equation}
\label{hatpat} \partial_t  f(t,\omega):=
\lim_{h\to0, h>0} \frac
{1}{h} \bigl[f (t+h,\omega^t )-f (t, \omega) \bigr],\qquad t<T,
\end{equation}
exists. If $ \partial_t  f(t,\omega)$ exists for all  $(t,\omega)\in\Lambda$, then  the non-anticipating functional $\partial f_t$ is called the
horizontal derivative of $f$.

A  non-anticipative functional $f(t,\omega)$ is vertically differentiable at $(t,\omega)\in\Lambda$ if
\begin{equation}
\label{hatpax} \partial_{\omega} f(t,\omega):=
\lim_{h\to0}\frac
{1}{h} \bigl[  f(t,\omega+ h
1_{[t,T]}) -  f(t,\omega) \bigr],
\end{equation}
exists. If $f$ is vertically differentiable at all $(t,\omega)\in\Lambda$ then the map  $\partial _\omega f :\L\to R$ defines a non-anticipative map,
 called the vertical derivative of $f$.

Similarly one can define
\begin{equation}
\partial_{\omega\omega}f:= \partial_{\omega
}(\partial f_{\omega}),\qquad.
\end{equation}

Define $C^{1,k}([0, T )\times \Omega)$ as the set of functionals $f$, which
are
\begin{itemize}
  \item  horizontally differentiable with $\partial_t f$ continuous at fixed times,
  \item $k$ times vertically differentiable with continuous $\partial_{\o}^k f$.
\end{itemize}
The following assertion follows from the generalized It\^o formula for cadlag semimartingales proved in  \cite{CF1} (see also  \cite{LScS}).

\be{thr}
Let $f\in C^{1,1}([0,T]\times \O)$. Then  for all $(t,\o)\in [0,T]\times \cal V$
\beaa
f(t,\o)=f(0,\o)+\int_0^t\partial_t f(s,\o)ds+\int_0^t\partial_\o f(s-,\o)d\o_s\\
+\sum_{s\le t}(f(s,\o)-f(s-,\o)-\partial_{\omega}f(s-,\o)\Delta\o_s)
\eeaa
and $f(t,\o)\in\cal V$ for all $\o\in\cal V$.
\ee{thr}
{\bf Corollary}. If $f\in C^{1,1}([0,T]\times \O)$, then $f$ is differentiable in the sense of Definition 1 and
$$
\partial_tf=f^0,\;\;\;\;\partial _\omega f= f^1.
$$

In order to compare Dupire's derivatives with Chitashvili's  derivative (the derivative in the sense of Definition 1), we introduce another type of
 vertical  derivative where, unlike to Dupire's derivative $\partial_\omega f$,  the path deformation  of
continuous paths remain continuous.

{\bf Definition 3}.
We say that a  non-anticipative functional $f(t,\omega)$ is  vertically differentiable   and denote this differential by $D_\omega f(t,\omega)$, if
the limit
\begin{equation}
D_\omega f(t,\omega):=\lim_{h\to0, h>0}\frac{f(t+h,\o^{t}+\chi_{t,h})-f(t+h,\o^{t})}{h},
\end{equation}
exists  for all $(t,\o)\in [0,T]\times \Omega$, where
$$
\chi_{t,h}(s)=(s-t)1_{(t,t+h]}(s)+h1_{(t+h,T]}(s).
$$
The second order derivative is defined similarly
$$
D_{\o\o}f=D_\o(D_\o f).
$$

Note that, if $f(t,\o)=g(\o_t)$ for any $\o\in D[0,T]$, where $g=(g(x), x\in R)$ is a differentiable function, then
$D_\o f(t,\o)$ (so as $\partial _\o f(t,\o)$) coincides with $g'(\o_t)$.

\be{prop}\lbl{11}
Let  $f\in C([0,T]\times \O)$  be  differentiable in the sense of Definition 1, i.e.,  there exist $f^0, f^1\in C([0,T]\times \O)$,
such that for all  $(t,\o)\in [0,T]\times\cal V$
\begin{equation}\label{xv22}
f(t,\o)=f(0,\o)+\int_0^tf^0(s,\o)ds+\int_0^tf^1(s-,\o)d\o_s + V(t,\o),
\end{equation}
where
$$
V(t,\o):=\sum_{s\le t}\big[f(s,\o)-f(s-,\o)-f^1(s-,\o)\Delta\o_s\big]
$$
is of finite variation for all  $\o\in {\cal V}$.

Then for all  $(t,\o)\in [0,T]\times D([0,T])$
\begin{equation}
f^0(t,\o)=\partial_t f(t,\omega)\;\;\;\;\text{and}\;\;\;\; f^1(t,\o)=D_\omega f(t,\omega).
\end{equation}
\ee{prop}

{\it Proof.}  Since $\o^t$ is constant on $[t,T]$ and $f(t,\o^t)=f(t,\o)$, if $s\le t$, from (\ref{xv22}) we have that for any $ \o\in{\cal V}$

\begin{equation}\label{xv23}
f(t+h,\o^t)=f(0,\o)+\int_0^tf^0(s,\o)ds+\int_t^{t+h}f^0(s,\o^t)ds+
\end{equation}
$$
+ \int_0^tf^1(s-,\o)d\o_s + V(t,\o)
$$
and
\begin{equation}\label{xv24}
f(t+h,\o^t+\chi_{t,h})=f(0,\o)+\int_0^tf^0(s,\o)ds+ \int_0^tf^1(s-,\o)d\o_s+
\end{equation}
$$
+\int_t^{t+h}f^0(s,\o^t+\chi_{t,h})ds+ \int_t^{t+h}f^1(s,\o^t+\chi_{t,h})ds+ V(t,\o).
$$
Therefore
$$
\partial_t f(t,\o)=\lim_{h\to0}\frac{f(t+h,\o^{t})-f(t,\o)}{h}=
$$
$$
=\lim_{h\to0}\frac{1}{h}\int_{t}^{t+h}f^0(s,\o^{t})ds= f^0(t,\o)
$$
by continuity of $f^0$.

It is evident that $\chi_{t,h}(s)\le h$  and
$$\frac{\chi_{t,h}(s)-\chi_{t,0}(s)}{h}=\frac{\chi_{t,h}(s)}{h}\to 1_{[t,T]}(s)\;as\;h\to0+,\;\forall s\in [0,T].$$

Trerefore, relations (\ref{xv24})-(\ref{xv23}) and continuity of $f^1$ and $f^0$ imply that
\beaa
D_\o f(t,\o)=\lim_{h\to0}\frac{f(t+h,\o^{t}+\chi_{t,h})-f(t+h,\o^{t})}{h}=\\
=\lim_{h\to0}\frac{1}{h}\int_{t}^{t+h}\big(f^0(s,\o^{t}+\chi_{t,h})-f^0(s,\o^{t})\big)ds\\
+\lim_{h\to0}\frac{1}{h}\int_{t}^{t+h}f^1(s,\o^{t}+\chi_{t,h})ds= f^1(t,\o)
\eeaa
for any $\o\in\cal V([0,T])$  and by continuity of $f^1$ this equality is true for all $\o\in D([0,T])$.

{\bf Remark.}
If  $f\in C([0,T]\times \O)$  is two times  differentiable in the sense of Definition 1, then similarly one can show that
$$
f^{1,1}(t,\o)=D_{\o\o}f(s,\o).
$$

{\bf Corrolary 1.} Let $f\in C^{1,1}([0,T]\times \O)$. Then for all $(t,\o)\in \Lambda$
\beaa
\partial_\o f(t,\o)=f^1(t,\o)=D_\o f(t,\o).
\eeaa

In general  $ \partial_\o f(t,\o)$ and  $D_\o f(t,\o)$  are not equal.

{\bf Counterexample 1}. Let $g=(g(x),x\in r)$  be a bounded differentiable function and let  $f(t,\o)=g(\o_t)-g(\o_{t-})$. Then $\partial_\o f(t,\o)=g'(\o_t)$
and
\beaa
D_\o f(t,\o)=\lim_{h\to 0+}\frac{f(t+h,\o^{t}+\chi_{t,h})-f(t+h,\o^{t})}{h}=0,\;
\rm{since}
\\f(t+h,\o^{t}+\chi_{t,h})-f(t+h,\o^{t})=g(\o_t+h)-g(\o_t+h)-g(\o_t)+g(\o_t)=0.
\eeaa
It is evident that $f\bar\in C^{1,1}(\L)$, since $f\bar\in C(\L)$ and
$\partial_t f=\infty.$

The following assertion shows that if $f$ belongs to the class $C^{1,2}(\Lambda)$ of non-anticipative functionals, then $\partial f_\o (t,\o)$ and  $\partial f_{\o\o} (t,\o)$
 are uniquelly determined  by the restriction of $f$ to continuous paths.  This assertion  is proved  by Cont and Fournie \cite{CF} (see also \cite{BCC})  in a complicated way.
We give a simple proof based on Proposition 1.

{\bf Corrolary 2.} Let $f^1$ and $f^2$ belong to $\in C^{1,2}(\Lambda)$ in the Dupire sense and
\begin{equation}\label{f1f2}
f^1(t,\o)=f^2(t,\o)\;\;\;\;\text{for all}\;\;\;(t,\o)\in [0,T]\times C([0,T]).
\end{equation}
Then
\begin{equation}\label{f12}
\partial_\o f^1(t,\o)=\partial_\o f^2(t,\o),\;\;\;\partial_{\o\o} f^1(t,\o)=\partial_{\o\o}f^2(t,\o)
\end{equation}
for all $(\o,t)\in [0,T]\times C([0,T])$.

{\it Proof}. By Theorem 2  
\begin{equation}\label{bv}
f^i(t,\o)=f^i(0,\o)+\int_0^t\partial_t f^i(s,\o)ds+\int_0^t\partial_\o f^i(s,\o)d\o_s\;\;\;i=1,2,
\end{equation}
 for all $\o\in C([0,T])\cap{\cal V}([0,T])$.
It follows from Proposition 1 that
$$
\partial_\o f^i(t,\o)=D_\o f^i(t,\o);\;\;\; i=1,2.
$$
 Since  $\o^{t}+\chi_{t,h}\in C([0,T])$ if $\o\in C([0,T])$, by  definition of $D_\o$ and equality (\ref{f1f2}) we have
\begin{equation}\label{df12}
D_\o f^1(t,\o)=D_\o f^2(t,\o)\;\;\;\;\text{for all}\;\;\;(t,\o)\in [0,T]\times C([0,T])),
\end{equation}
which implies that
\begin{equation}\label{f12}
\partial_\o f^1(t,\o)=\partial_\o f^2(t,\o),\;\;\;\;\text{for all}\;\;\;(t,\o)\in [0,T]\times C([0,T]),
\end{equation}
It is evident that  $\partial_t f^1(t,\o)=\partial_t f^2(t,\o)$ for all $(\o,t)\in [0,T]\times C([0,T])$. Therefore, comparing the It\^o formulas (\ref{itoc}) for
$f^1(t,\o)$ and $f^2(t,\o)$ we obtain that
$$
\int_t^u\partial_{\o\o}f^1(s,\o)d\la\o\ra_s=\int_t^u\partial_{\o\o}f^2(s,\o)d\la\o\ra_s
$$
for any continuous semimartinale $\o$. Dividing both parts of this equality by $\la\o\ra_u-\la\o\ra_t$ and passing to the limit as $u\to t$, we obtain that
$\partial_{\o\o} f^1(t,\o)=\partial_{\o\o}f^2(t,\o)$ for any continuous semimartingale and by continuity of $\partial_{\o\o} f^1(t,\o)$ and $\partial_{\o\o}f^2(t,\o)$ this equality
will be true for all $\o\in C([0,T])$.

\be{prop}\lbl{33}
Let $f\in C([0,T]\times \O)$ be differentiable in the sense of Definition 1 and
\begin{equation}\label{xv25}
\left|f(t,\o)-f(t-,\o)-\Delta\o_tf^1(t-,\o)\right|\le K|\Delta\o_t|^2
\end{equation}
for some $K>0$.
Then
$$
f^0(t,\o)=\partial_t f(t,\o), \;\;\;\;\forall(t,\o)\in \L,
$$
$$
f^1(t,\o)=\partial_\o f(t,\o),\;\;\;\;\forall \o\in C[0,T]
$$
(or for all $\omega$ continuous at $t$).
\ee{prop}

{\it Proof.} For $\o\in D[0,T]$ let
$\tilde\o=\o_s$ if $s<t$ and $\tilde\o=\o_{s-}+h$, if $s\ge t$, i.e. $\td\o=\o^{t-}+h1_{[t,T]}$, hence $\Delta\td\o_s=h$.

Therefore, using condition (\ref{xv25}) for $\tilde\o$ we have
\beaa
\left|\frac{f(t,\o^{t-}+h1_{[t,T]})-f(t-,\o)}{h}-f^1(t-,\o)\right|\le K|h|,\;\forall h.
\eeaa
It follows from here that
$$
\lim_{h\to0}\frac{f(t,\o^{t-}+h1_{[t,T]})-f(t-,\o)}{h}=f^1(t-,\o),
$$
which implies that
$f^1(t,\o)=\partial_\o f(t,\o)$ if $\o$ is continuous at $t$.
Equality $f^1(t,\o)=\partial_tf(t,\o), \forall(t,\o)\in \L$ is proved in Proposition 1.\qed

Now we introduce definition of space derivatives  which can be calculated pathwise along the differentiable paths
and   using such derivatives in  Theorem 3  below  a change of variables formula for functions of finite variations  is proved, which
gives sufficient conditions for the existence of derivatives in the Chitashvili sense.

{\bf Definition 4}.
We say that a  non-anticipative functional $f(t,\omega)$ is  differentiable, if the limits  $f_t\;f_\o\in C(\L)$ exist, where
\beaa
f_t(t,\o)=\lim_{h\to0, h>0}\frac{f(t+h,\o^{t})-f(t,\o)}{h}, \;\;\;\; \forall(t,\o)\in [0,T]\times D[0,T]\\
f_{\o}(t,\o)=\lim_{h\to0, h>0}\frac{f(t+h,\o)-f(t+h,\o^{t})}{\o_{t+h}-\o_t},\;\;\;\;\forall(t,\o)\in [0,T]\times C^1[0,T].\\
\eeaa

\be{prop}\lbl{22}
Let $f$ be differentiable in the sense of definition 4.
Then $\forall(t,\o)\in [0,T]\times C^1[0,T]$
\beq\lbl{itt}
f(t,\o)-f(0,\o)=\int_0^tf_t(s,\o)ds+\int_0^tf_\o(s,\o)d\o_s.
\eeq
\ee{prop}
{\it Proof}. We have
\beaa
\lim_{h\to0, h>0}\frac{f(t+h,\o)-f(t,\o)}{h}\\
=
\lim_{h\to0+}\frac{f(t+h,\o)-f(t+h,\o^{t})}{{\o_{t+h}-\o_t}}\times \frac{{\o_{t+h}-\o_t}}{h}
\\
+\lim_{h\to0+}\frac{f(t+h,\o^{t})-f(t,\o)}{h}
=\o'(t)f_\o(t,\o)+f_t(t,\o),\\
\forall(t,\o)\in [0,T]\times C^1[0,T].
\eeaa
Hence right derivative of
\beaa
f(t,\o)-f(0,\o)-\int_0^tf_t(s,\o)ds-\int_0^tf_\o(s,\o)\o'_sds
\eeaa
is zero for each $\o\in C^1$. By the Lemma A2 of appendix formula (\ref{itt}) is satisfied.

\be{thr}
Let $f\in C(\L)$ and $f_t,f_\o\in C(\L)$ are derivatives in the sense of definition 4.
Assume  also that  for any $\o\in\cal V$
$$
\sum_{s\le t}|f(s,\o)-f(s-,\o)|<\infty.
$$
Then
\beaa
f(t,\o)=f(0,\o)+\int_0^tf_t(s,\o)ds+\int_0^tf_\o(s,\o)d\o_s^c\\
+\sum_{s\le t}(f(s,\o)-f(s-,\o)).
\eeaa
\ee{thr}
{\it Proof}.
For $\o\in V$ we have $\o=\o^c+\o^d,\;\o^d=\sum_{s\le t}\Delta\o_s,\;\o^c\in C$. Set
$$\o^{d,n}=\sum_{s\le t,|\Delta\o_s|>\frac1n}\Delta\o_s,\;\o^n=\o^c+\o^{d,n}.$$
It is evident that as $n\to 0$
$$
|\o^n-\o|_T=|\o^d-\o^{d,n}|_T=\max_t|\int_0^t1_{(|\Delta\o_s|\le\frac1n)}d\o_s^d\le\int_0^T1_{(|\Delta\o_s|\le\frac1n)}dvar_s(\o^d)\to 0.
$$
 We know that discontinuity points of $f$ are also discontinuity points of $\o$.
Let $\{t_1<...<t_k\}=\{s:|\Delta\o_s|>\frac1n\}\cup\{0,T\}$. Denote by $\o^{\e}\in C'$ a differentiable approximation of $\o^c$,  such that $var_T(\o^\e-\o^c)<\e$  and
let $\o^{n,\e}=\o^\e+\o^{d,n}$.
Then by  Proposition \ref{22}
\beaa
f(t,\o^{n,\e})-f(t_i,\o^{n,\e})-\int_{t_i}^tf_t(s,\o^{n,\e})ds-\int_{t_i}^tf_\o^{n,\e}(s,\o^{n,\e})\o^{'\e}_sds=0,\;t\in[t_i,t_{i+1})
\eeaa
and
\beaa
f(T,\o^{n,\e})-f(0,\o^{n,\e})=\sum_{i\ge 1} \big(f(t_{i},\o^{n,\e})-f(t_{i-1},\o^{n,\e})\big)\\
=\sum \big (f(t_{i}-,\o^{n,\e})-f(t_{i-1},\o^{n,\e})\big )+\sum \big (f(t_i,\o^{n,\e})-f(t_i-,\o^{n,\e})\big )\\
=\sum\int_{t_{i-1}}^{t_i}f_t(s,\o^{n,\e})ds+\sum\int_{t_{i-1}}^{t_i}f_\o(s,\o^{n,\e})\o^{'\e}_sds+\sum \big (f(t_i,\o^{n,\e})-f(t_i-,\o^{n,\e})\big )\\
=\sum\int_{t_{i-1}}^{t_i}f_t(s,\o^{n,\e})ds+\sum\int_{t_{i-1}}^{t_i}f_\o(s,\o^{n,\e})d\o_s^{n,\e}\\
-\sum f_\o(t_i-,\o^{n,\e})\Delta\o^{n,\e}_{t_i}+\sum \big (f(t_i,\o^{n,\e})-f(t_i-,\o^{n,\e})\big )\\
=\int_0^Tf_s(s,\o^{n,\e})ds+\int_0^Tf_\o(s,\o^{n,\e})d\o^{n,\e}_s\\
+\sum\big ( f(t_i,\o^{n,\e})-f(t_i-,\o^{n,\e})-f_\o(t_i-,\o^{n,\e})\Delta\o^{n,\e}_{t_i}\big )\\
=\int_0^Tf_t(s,\o^{n,\e})ds+\int_0^Tf_\o(s,\o^{n,\e})d\o_s^\e+\sum\big ( f(t_i,\o^{n,\e})-f(t_i-,\o^{n,\e})\big ).
\eeaa
Since $f(t,\o^{n,\varepsilon})$  admits finite number of jumps and $\sup_\varepsilon var_T\o^\varepsilon<\infty$, passing to the limit as $\e\to0$ we get
\beaa
f(T,\o^{n})-f(0,\o^{n})\\
=\int_0^Tf_t(s,\o^{n})ds+\int_0^Tf_\o(s,\o^{n})d\o_s^c+\sum\big ( f(t_i,\o^n)-f(t_i-,\o^n)\big ).
\eeaa
By the continuity of functionals
$f,\;f_t,\;f_\o$ and Lemma A1 from the appendix
$$f(t,\o^n)\to f(t,\o),\;\int_0^tf_t(s,\o^n)ds\to \int_0^tf_t(s,\o)ds,
$$
$$\int_0^tf_t(s,\o^n)d\o_s^c\to \int_0^tf_t(s,\o)d\o_s^c,\;as\;n\to\infty.$$
 It remains to show   convergence of the sum. Since\\
 $f^d(t,\o)=\sum_{s\le t}f(s,\o)-f(s-,\o)$ is  of finite variation
 \beaa
f^d(t,\o)=\sum \big ( f(t_i,\o^{n})-f(t_i-,\o^{n})\big )-\sum \big (f(t_i,\o)-f(t_i-,\o)\big )\\
=\sum_{s\le t} (f(s,\o)-f(s-,\o))1_{(|\Delta\o_s|\le\frac1n)}\\
=\int_0^t1_{(|\Delta\o_s|\le\frac1n)}df^d(s,\o)\to 0,\; as\;n\to\infty,
 \eeaa
by  the dominated convergence theorem.

{\bf Corollary.} If $f$ satisfies conditions of Theorem 3 then $f$ is differentiable in the sense of Definition 1.

\section{Appendix}

The following  lemma is a modification of lemma 6 of \cite{LScS}.

{\bf Lemma A1}. Let $X_n,X\in \O$ be a  sequence of paths, such that $||X_n-X||_T\to 0$ as $n\to\infty$.   Let $f\in C(\L)$.
Then
$$\sup_{t\le T}|f(t,X_n)-f(t,X)|\overset{n\to\infty}\to 0.$$

{\it Proof}.
If not then $\exists \e > 0$, a sequence of integers $n_k, k=1,...$, and a sequence $s_k\in [0, T ]$ such that
\beq\lbl{uni}
|f(s_k, X_{n_k})-f(s_k, X)|
\ge \e
\eeq
By moving to a subsequence we can assume without loss of generality that either $s_k \to s^*,\;s_k\ge s^*$  or $s_k \to s^*,\;s_k< s^*$ for some $s^*\in[0, T]$.
In the first case by continuity assumption we get
\beaa
|f(s_k, X_{n_k})-f(s_k, X)|\le |f(s_k, X_{n_k})-f(s^*, X)|\\+|f(s_k, X)-f(s^*, X)|
\to 0,
\eeaa
since $d_\infty((s_k,X_{n_k}),(s^*, X))\to 0,\;d_\infty((s_k,X),(s^*, X))\to 0$ .

In the second case we have
\beaa
|f(s_k, X_{n_k})-f(s_k, X)|
\le |f(s_k, X_{n_k})-f(s^*, X^{s^*-})|\\
+|f(s_k, X)-f(s^*, X^{s^*-})|\to 0,
\eeaa
since $d_\infty((s_k, X_{n_k}),(s^*, X^{s^*-}))\to 0,\;d_\infty((s_k, X),(s^*, X^{s^*-}))\to 0$ .
This contradicts (\ref{uni}).

We shall need also the following assertion

{\bf Lemma A2}. Let $f$ be a real-valued, continuous function, defined on an arbitrary interval $I$ of the real line. If $f$ is right (or left) differentiable
at every point $a \in I$, which is not the supremum (infimum) of the interval, and if this right (left) derivative is always zero, then $f$ is a constant.

{\it Proof}.
For a proof by contradiction, assume there exist $a < b$ in $I$ such that $f(a) \neq f(b)$. Then
\beaa
   \varepsilon :={\frac {|f(b)-f(a)|}{2(b-a)}}>0.
\eeaa
Define $c$ as the infimum of all those $x$ in the interval $(a,b]$ for which the difference quotient of $f$ exceeds $\varepsilon$ in absolute value, i.e.
\beaa
     c=\inf\{\,x\in (a,b]\mid |f(x)-f(a)|>\varepsilon (x-a)\,\}.
\eeaa
Due to the continuity of $f$, it follows that $c < b$ and $|f(c)-–f(a)|=\varepsilon(c–-a)$. At $c$ the right derivative of $f$ is zero by assumption, hence there exists $d$ in the interval $(c,b]$
with $|f(x)–-f(c)|\le\varepsilon(x–-c)$ for all $x \in (c,d]$. Hence, by the triangle inequality,
\beaa
      |f(x)-f(a)|\leq |f(x)-f(c)|+|f(c)-f(a)|\leq \varepsilon (x-a)
\eeaa
for all $x$ in $[c,d)$, which contradicts the definition of $c$.

\end{document}